\documentclass{article}%
\usepackage{amsmath}
\usepackage{amsfonts}
\usepackage{amssymb}
\usepackage{graphicx}%
\providecommand{\U}[1]{\protect \rule{.1in}{.1in}}
\newtheorem{theorem}{Theorem}

\newtheorem{definition}[theorem]{Definition}
\newtheorem{example}[theorem]{Example}

\newtheorem{lemma}[theorem]{Lemma}

\newtheorem{proposition}[theorem]{Proposition}
\newtheorem{remark}[theorem]{Remark}

\newenvironment{proof}[1][Proof]{\noindent \textbf{#1.} }{\  \rule{0.5em}{0.5em}}
\begin{document}

\title{$G$--Expectation, $G$--Brownian Motion and Related Stochastic Calculus of
It\^{o} Type}
\author{Shige PENG\thanks{The author thanks the partial support from the Natural
Science Foundation of China, grant No. 10131040. He thanks to the anonymous
referee's constructive suggestions, as well as Juan Li's typos-corrections.
Special thanks are to the organizers of the memorable Abel Symposium 2005 for
their warm hospitality and excellent work (see also this paper in:
http://abelsymposium.no/2005/preprints). }\\Institute of Mathematics, Fudan University\\Institute of Mathematics\\Shandong University\\250100, Jinan, China\\peng@sdu.edu.cn}
\date{1st version: arXiv:math.PR/0601035 v1 3 Jan 2006}
\maketitle

\noindent \  \noindent \textbf{Abstract. }{\small We introduce a notion of
nonlinear expectation ----}$G${\small --expectation---- generated by a
nonlinear heat equation with a given infinitesimal generator }$G${\small . We
first discuss the notion of }$G${\small --standard normal distribution. With
this nonlinear distribution we can introduce our }$G${\small --expectation
under which the canonical process is a }$G${\small --Brownian motion. We then
establish the related stochastic calculus, especially stochastic integrals of
It\^{o}'s type with respect to our }$G${\small --Brownian motion and derive
the related It\^{o}'s formula. We have also given the existence and uniqueness
of stochastic differential equation under our }$G${\small --expectation. As
compared with our previous framework of }$g${\small --expectations, the theory
of }$G${\small --expectation is intrinsic in the sense that it is not based on
a given (linear) probability space.}\newline \newline


\noindent \textbf{MSC 2000 Classification Numbers: }60H10, 60H05, 60H30, 60J60,
60J65, 60A05, 60E05, 60G05, 60G51, 35K55, 35K15, 49L25\newline \newline

\section{Introduction}

In 1933 Andrei Kolmogorov published his Foundation of Probability Theory
(Grundbegriffe der Wahrscheinlichkeitsrechnung) which set out the axiomatic
basis for modern probability theory. The whole theory is built on the Measure
Theory created by \'{E}mile Borel and Henry Lebesgue and profoundly developed
by Radon and Fr\'{e}chet. The triple $(\Omega,\mathcal{F},\mathbf{P})$, i.e.,
a measurable space $(\Omega,\mathcal{F)}$ equipped with a probability measure
$\mathbf{P}$ becomes a standard notion which appears in most papers of
probability and mathematical finance. The second important notion, which is in
fact at an equivalent place as the probability measure itself, is the notion
of expectation. The expectation $\mathbf{E}[X]$ of a $\mathcal{F}$--measurable
random variable $X$ is defined as the integral $\int_{\Omega}XdP$. A very
original idea of Kolmogorov's Grundbegriffe is to use Radon--Nikodym theorem
to introduce the conditional probability and the related conditional
expectation under a given $\sigma$--algebra $\mathcal{G\subset F}$. It is hard
to imagine the present state of arts of probability theory, especially of
stochastic processes, e.g., martingale theory, without such notion of
conditional expectations. A given time information $(\mathcal{F}_{t})_{t\geq
0}$ is so ingeniously and consistently combined with the related conditional
expectations $\mathbf{E}[X|\mathcal{F}_{t}]_{t\geq0}$. It\^{o}'s
calculus---It\^{o}'s integration, It\^{o}'s formula and It\^{o}'s equation
since 1942 \cite{Ito1942}, is, I think, the most beautiful discovery on this ground.

A very interesting problem is to develop a nonlinear expectation
$\mathbb{E}[\cdot]$ under which we still have such notion of conditional
expectation. A notion of $g$--expectation was introduced by Peng, 1997
(\cite{Peng1997} and \cite{Peng1997b}) in which the conditional expectation
$\mathbb{E}^{g}[X|\mathcal{F}_{t}]_{t\geq0}$ is the solution of the backward
stochastic differential equation (BSDE), within the classical framework of
It\^{o}'s calculus, with $X$ as its given terminal condition and with a given
real function $g$ as the generator of the BSDE. driven by a Brownian motion
defined on a given probability space $(\Omega,\mathcal{F},\mathbf{P})$. It is
completely and perfectly characterized by the function $g$. The above
conditional expectation is characterized by the following well-known
condition.
\[
\mathbb{E}^{g}[\mathbb{E}^{g}[X|\mathcal{F}_{t}]\mathbf{I}_{A}]=\mathbb{E}%
^{g}[X\mathbf{I}_{A}],\  \  \  \forall A\in \mathcal{F}_{t}.
\]
Since then many results have been obtained in this subject (see, among others,
\cite{BCHMP1}, \cite{Chen98}, \cite{CE}, \cite{CKJ}, \cite{CHMP},
\cite{CHMP3}, \cite{CP}, \cite{CP1}, \cite{Jiang}, \cite{JC}, \cite{Peng1999},
\cite{Peng2005a}, \cite{Peng2005b}, \cite{PX2003}, \cite{Roazza2003},
\cite{K-S}).

In \cite{Peng2005} (see also \cite{Peng2004}), we have constructed a kind of
filtration--consistent nonlinear expectations through the so--called nonlinear
Markov chain. As compared with the framework of $g$--expectation, the theory
of $G$--expectation is intrinsic, a meaning similar to the \textquotedblleft
intrinsic geometry\textquotedblright. in the sense that it is not based on a
classical probability space given a priori.

In this paper, we concentrate ourselves to a concrete case of the above
situation and introduce a notion of $G$--expectation which is generated by a
very simple one dimensional fully nonlinear heat equation, called $G$--heat
equation, whose coefficient has only one parameter more than the classical
heat equation considered since Bachelier 1900, Einstein 1905 to describe the
Brownian motion.. But this slight generalization changes the whole things.
Firstly, a random variable $X$ with \textquotedblleft$G$--normal
distribution\textquotedblright \ is defined via the heat equation. With this
single nonlinear distribution we manage to introduce our $G$--expectation
under which the canonical process is a $G$--Brownian motion.

We then establish the related stochastic calculus, especially stochastic
integrals of It\^{o}'s type with respect to our $G$--Brownian motion. A new
type of It\^{o}'s formula is obtained. We have also established the existence
and uniqueness of stochastic differential equation under our $G$--stochastic calculus.

In this paper we concentrate ourselves to $1$--dimensional $G$--Brownian
motion. But our method of \cite{Peng2005} can be applied to multi--dimensional
$G$--normal distribution, $G$--Brownian motion and the related stochastic
calculus. This will be given in \cite{Peng2006b}.

Recently a new type of second order BSDE was proposed to give a probabilistic
approach for fully nonlinear 2nd order PDE, see \cite{Touzi}. In finance a
type of uncertain volatility model in which the PDE of Black-Scholes type was
modified to a fully nonlinear model, see \cite{Avellaneda} and \cite{Lyons}. A
point of view of nonlinear expectation and conditional expectation was
proposed in \cite{Peng2004} and \cite{Peng2005}. When I presented the result
of this paper in Workshop on Risk Measures in Evry, July 2006, I met Laurent
Denis and got to learn his interesting work, joint with Martini, on volatility
model uncertainty \cite{Denis}. See also our forthcoming paper
\cite{Denis-Peng} for the pathwise analysis of $G$-Brownian motion.

As indicated in Remark \ref{Rem-1}, the nonlinear expectations discussed in
this paper are equivalent to the notion of coherent risk measures. This with
the related conditional expectations $\mathbb{E}[\cdot|\mathcal{F}_{t}%
]_{t\geq0}$ makes a dynamic risk measure: $G$--risk measure.

This paper is organized as follows: in Section 2, we recall the framework
established in \cite{Peng2005} and adapt it to our objective. In section 3 we
introduce $1$--dimensional standard $G$-normal distribution and discuss its
main properties. In Section 4 we introduce 1--dimensional $G$-Brownian motion,
the corresponding $G$--expectation and their main properties. We then can
establish stochastic integral with respect to our $G$-Brownian motion of
It\^{o}'s type and the corresponding It\^{o}'s formula in Section 5 and the
existence and uniqueness theorem of SDE driven by $G$-Brownian motion in
Section 6.

\section{Nonlinear expectation: a general framework}

We briefly recall the notion of nonlinear expectations introduced in
\cite{Peng2005}. Following Daniell (see Daniell 1918 \cite{Daniell}) in his
famous Daniell's integration, we begin with a vector lattice. Let $\Omega$ be
a given set and let $\mathcal{H}$ be a vector lattice of real functions
defined on $\Omega$ containing $1$, namely, $\mathcal{H}$ is a linear space
such that $1\in \mathcal{H}$ and that $X\in \mathcal{H}$ implies $|X|\in
\mathcal{H}$. $\mathcal{H}$ is a space of random variables. We assume the
functions on $\mathcal{H}$ are all bounded. Notice that
\[
a\wedge b=\min \{a,b\}=\frac{1}{2}(a+b-|a-b|),\  \ a\vee b=-[(-a)\wedge(-b)].
\]
Thus $X$, $Y\in \mathcal{H}$ implies that $X\wedge Y$, $X\vee Y$, $X^{+}%
=X\vee0$ and $X^{-}=(-X)^{+}$ are all in $\mathcal{H}$.

\begin{definition}
\label{Def-1}\textbf{A nonlinear expectation }$\mathbb{E}$ is a functional
$\mathcal{H}\mapsto \mathbb{R}$ satisfying the following properties\newline%
\newline \textbf{(a) Monotonicity:} If $X,Y\in \mathcal{H}$ and $X\geq Y$ then
$\mathbb{E}[X]\geq \mathbb{E}[Y].$\newline \textbf{(b)} \textbf{Preserving of
constants:} $\mathbb{E}[c]=c$.\newline \newline In this paper we are interested
in the expectations which satisfy\newline \newline \textbf{(c)}
\textbf{Sub-additivity (or self--dominated property):}%
\[
\mathbb{E}[X]-\mathbb{E}[Y]\leq \mathbb{E}[X-Y],\  \  \forall X,Y\in \mathcal{H}.
\]
\textbf{(d) Positive homogeneity: } $\mathbb{E}[\lambda X]=\lambda
\mathbb{E}[X]$,$\  \  \forall \lambda \geq0$, $X\in \mathcal{H}$.\newline%
\textbf{(e) Constant translatability: }$\mathbb{E}[X+c]=\mathbb{E}[X]+c$.
\end{definition}

\medskip

\begin{remark}
\label{Rem-Def1}The above condition (d) has an equivalent form: $\mathbb{E}%
[\lambda X]=\lambda^{+}\mathbb{E}[X]+\lambda^{-}\mathbb{E}[-X]$. This form
will be very convenient for the conditional expectations studied in this paper
(see (vi) of Proposition \ref{Prop-1-7}).
\end{remark}

\begin{remark}
\label{Rem-1}We recall the notion of the above expectations satisfying
(c)--(e) was systematically introduced by Artzner, Delbaen, Eber and Heath
\cite{ADEH1}, \cite{ADEH2}, in the case where $\Omega$ is a finite set, and by
Delbaen \cite{Delbaen} in general situation with the notation of risk measure:
$\rho(X)=\mathbb{E}[-X]$. See also in Huber \cite{Huber} for even early study
of this notion $\mathbb{E}$ (called upper expectation $\mathbf{E}^{\ast}$ in
Ch.10 of \cite{Huber}) in a finite set $\Omega$. See Rosazza Gianin
\cite{Roazza2003} or Peng \cite{Peng2003}, El Karoui \& Barrieu \cite{El-Bar},
\cite{El-Bar2005} for dynamic risk measures using $g$--expectations.
Super-hedging and super pricing (see \cite{EQ} and \cite{EPQ}) are also
closely related to this formulation.
\end{remark}

\begin{remark}
\label{Rem-2}We observe that $\mathcal{H}_{0}=\{X\in \mathcal{H}$,
$\mathbb{E}[|X|]=0\}$ is a linear subspace of $\mathcal{H}$. To take
$\mathcal{H}_{0}$ as our null space, we introduce the quotient space
$\mathcal{H}/\mathcal{H}_{0}$. Observe that, for every $\mathbf{\{}X\}
\in \mathcal{H}/\mathcal{H}_{0}$ with a representation $X\in \mathcal{H}$, we
can define an expectation $\mathbb{E}\mathbf{[\{}X\}]:=\mathbb{E}[X]$ which
still satisfies (a)--(e) of Definition \ref{Def-1}. Following \cite{Peng2005},
we set $\left \Vert X\right \Vert :=\mathbb{E}[|X|]$, $X\in \mathcal{H}%
/\mathcal{H}_{0}$. It is easy to check that $\mathcal{H}/\mathcal{H}_{0}$ is a
normed space under $\left \Vert \cdot \right \Vert $. We then extend
$\mathcal{H}/\mathcal{H}_{0}$ to its completion $[\mathcal{H}]$ under\ this
norm. $([\mathcal{H}],\left \Vert \cdot \right \Vert )$ is a Banach space. The
nonlinear expectation $\mathbb{E}[\cdot]$ can be also continuously extended
from $\mathcal{H}/\mathcal{H}_{0}$ to $[\mathcal{H}]$, which satisfies (a)--(e).
\end{remark}

For any $X\in \mathcal{H}$, the mappings
\[
X^{+}(\omega):\mathcal{H\longmapsto H}\  \  \  \text{and \  \ }X^{-}%
(\omega):\mathcal{H\longmapsto H}%
\]
satisfy
\[
|X^{+}-Y^{+}|\leq|X-Y|\text{ \  \ and \ }\ |X^{-}-Y^{-}|=|(-X)^{+}%
-(-Y)^{+}|\leq|X-Y|.
\]
Thus they are both contraction mappings under $\left \Vert \cdot \right \Vert $
and can be continuously extended to the Banach space $([\mathcal{H}%
],\left \Vert \cdot \right \Vert )$.

We define the partial order \textquotedblleft$\geq$\textquotedblright \ in this
Banach space.

\begin{definition}
An element $X$ in $([\mathcal{H}],\left \Vert \cdot \right \Vert )$ is said to be
nonnegative, or $X\geq0$, $0\leq X$, if $X=X^{+}$. We also denote by $X\geq
Y$, or $Y\leq X$, if $X-Y\geq0$.
\end{definition}

It is easy to check that $X\geq Y$ and $Y\geq X$ implies $X=Y$ in
$([\mathcal{H}],\left \Vert \cdot \right \Vert )$.

The nonlinear expectation $\mathbb{E}[\cdot]$ can be continuously extended to
$([\mathcal{H}],\left \Vert \cdot \right \Vert )$ on which \textbf{(a)--(e)}
still hold.

\section{$G$--normal distributions}

For a given positive integer $n$, we denote by $lip(\mathbb{R}^{n})$ the space
of all bounded and Lipschitz real functions on $\mathbb{R}^{n}$. In this
section $\mathbb{R}$ is considered as $\Omega$ and $lip(\mathbb{R})$ as
$\mathcal{H}$.

In classical linear situation, a random variable $X(x)=x$ with standard normal
distribution, i.e., $X\thicksim N(0,1)$, can be characterized by
\[
E[\phi(X)]=\frac{1}{\sqrt{2\pi}}\int_{-\infty}^{\infty}e^{-\frac{x^{2}}{2}%
}\phi(x)dx,\  \  \forall \phi \in lip(\mathbb{R}).
\]
It is known since Bachelier 1900 and Einstein 1950 that $E[\phi(X)]=u(1,0)$
where $u=u(t,x)$ is the solution of the heat equation
\begin{equation}
\partial_{t}u=\frac{1}{2}\partial_{xx}^{2}u \label{heat0}%
\end{equation}
with Cauchy condition $u(0,x)=\phi(x)$.

In this paper we set $G(a)=\frac{1}{2}(a^{+}-\sigma_{0}^{2}a^{-})$,
$a\in \mathbb{R}$, where $\sigma_{0}\in \lbrack0,1]$ is fixed.

\begin{definition}
A real valued random variable $X$ with the standard $G$\textbf{--normal
distribution} is characterized by its $G$--expectation defined by
\[
\mathbb{E}[\phi(X)]=P_{1}^{G}(\phi):=u(1,0),\  \  \phi \in lip(\mathbb{R}%
)\mapsto \mathbb{R}%
\]
where $u=u(t,x)$ is a bounded continuous function on $[0,\infty)\times
\mathbb{R}$ which is the (unique) viscosity solution of the following
nonlinear parabolic partial differential equation (PDE)
\begin{equation}
\partial_{t}u-G(\partial_{xx}^{2}u)=0,\ u(0,x)=\phi(x). \label{eq-heat}%
\end{equation}
\end{definition}

\bigskip In case no confusion is caused, we often call the functional
$P_{1}^{G}(\cdot)$ the standard $G$--normal distribution. When $\sigma_{0}=1$,
the above PDE becomes the standard heat equation (\ref{heat0}) and thus this
$G$--distribution is just the classical normal distribution $N(0,1)$:%
\[
P_{1}^{G}(\phi)=P_{1}(\phi):=\frac{1}{\sqrt{2\pi}}\int_{-\infty}^{\infty
}e^{-\frac{x^{2}}{2}}\phi(x)dx.
\]

\begin{remark}
\label{Rem-3}The function $G$ can be written as $G(a)=\frac{1}{2}\sup
_{\sigma_{0}\leq \sigma \leq1}\sigma^{2}a$, thus the nonlinear heat equation
(\ref{eq-heat}) is a special kind of Hamilton--Jacobi--Bellman equation. The
existence and uniqueness of (\ref{eq-heat}) in the sense of viscosity solution
can be found in, for example, \cite{CIL}, \cite{FS}, \cite{Peng1992},
\cite{Yong-Zhou}, and \cite{Krylov} for $C^{1,2}$-solution if $\sigma_{0}>0$
(see also in \cite{Oksendal} for elliptic cases). Readers who are unfamililar
with the notion of viscosity solution of PDE can just consider, in the whole
paper, the case $\sigma_{0}>0$, under which the solution $u$ becomes a
classical smooth function.
\end{remark}

\begin{remark}
It is known that $u(t,\cdot)\in lip(\mathbb{R})$ (see e.g. \cite{Yong-Zhou}
Ch.4, Prop.3.1 or \cite{Peng1992} Lemma 3.1 for the Lipschitz continuity of
$u(t,\cdot)$, or Lemma 5.5 and Proposition 5.6 in \cite{Peng2004} for a more
general conclusion). The boundedness is simply from the comparison theorem (or
maximum principle) of this PDE. It is also easy to check that, for a given
$\psi \in lip(\mathbb{R}^{2})$, $P_{1}^{G}(\psi(x,\cdot))$ is still a bounded
and Lipschitz function in $x$.
\end{remark}

In general situations we have, from the comparison theorem of PDE,
\begin{equation}
P_{1}^{G}(\phi)\geq P_{1}(\phi),\  \forall \phi \in lip(\mathbb{R})\text{.}%
\  \label{compar}%
\end{equation}
The corresponding normal distribution with mean at $x\in \mathbb{R}$ and
variance $t>0$ is $P_{1}^{G}(\phi(x+\sqrt{t}\times \cdot))$. Just like the
classical situation, we have

\begin{lemma}
For each $\phi \in lip(\mathbb{R})$, the function
\begin{equation}
u(t,x)=P_{1}^{G}(\phi(x+\sqrt{t}\times \cdot)),\  \ (t,x)\in \lbrack
0,\infty)\times \mathbb{R} \label{u(t,x)}%
\end{equation}
is the solution of the nonlinear heat equation (\ref{eq-heat}) with the
initial condition $u(0,\cdot)=\phi(\cdot)$.
\end{lemma}

\begin{proof}
Let $u\in C([0,\infty)\times \mathbb{R})$ be the viscosity solution of
(\ref{eq-heat}) with $u(0,\cdot)=\phi(\cdot)\in lip(\mathbb{R})$. For a fixed
$(\bar{t},\bar{x})\in(0,\infty)\times \mathbb{R}$, we denote $\bar
{u}(t,x)=u(t\times \bar{t},x\sqrt{\bar{t}}+\bar{x})$. Then $\bar{u}$ is the
viscosity solution of (\ref{eq-heat}) with the initial condition $\bar
{u}(0,x)=\phi(x\sqrt{\bar{t}}+\bar{x})$. Indeed, let $\psi$ be a $C^{1,2}$
function on $(0,\infty)\times \mathbb{R}$ such that $\psi \geq \bar{u}$ (resp.
$\psi \leq \bar{u}$) and $\psi(\tau,\xi)=\bar{u}(\tau,\xi)$ for a fixed
$(\tau,\xi)\in(0,\infty)\times \mathbb{R}$. We have $\psi(\frac{t}{\bar{t}%
},\frac{x-\bar{x}}{\sqrt{\bar{t}}})\geq u(t,x)$, for all $(t,x)$ and
\[
\psi(\frac{t}{\bar{t}},\frac{x-\bar{x}}{\sqrt{\bar{t}}})=u(t,x),\  \text{at
}(t,x)=(\tau \bar{t},\xi \sqrt{\bar{t}}+\bar{x}).
\]
Since $u$ is the viscosity solution of (\ref{eq-heat}), at the point
$(t,x)=(\tau \bar{t},\xi \sqrt{\bar{t}}+\bar{x})$, we have
\[
\frac{\partial \psi(\frac{t}{\bar{t}},\frac{x-\bar{x}}{\sqrt{\bar{t}}}%
)}{\partial t}-G(\frac{\partial^{2}\psi(\frac{t}{\bar{t}},\frac{x-\bar{x}%
}{\sqrt{\bar{t}}})}{\partial x^{2}})\leq0\  \ (\text{resp. }\geq0).
\]
But since $G$ is positive homogenous, i.e., $G(\lambda a)=\lambda G(a)$, we
thus derive
\[
(\frac{\partial \psi(t,x)}{\partial t}-G(\frac{\partial^{2}\psi(t,x)}{\partial
x^{2}}))|_{(t,x)=(\tau,\xi)}\leq0\  \ (\text{resp. }\geq0).
\]
This implies that $\bar{u}$ is the viscosity subsolution (resp. supersolution)
of (\ref{eq-heat}). According to the definition of $P^{G}(\cdot)$ we obtain
(\ref{u(t,x)}).
\end{proof}

\begin{definition}
\label{Def-2}We denote
\begin{equation}
P_{t}^{G}(\phi)(x)=P_{1}^{G}(\phi(x+\sqrt{t}\times \cdot))=u(t,x),\  \ (t,x)\in
\lbrack0,\infty)\times \mathbb{R}. \label{P_tx}%
\end{equation}
\end{definition}

From the above lemma, for each $\phi \in lip(\mathbb{R})$, we have the
following Kolmogorov--Chapman chain rule:%
\begin{equation}
P_{t}^{G}(P_{s}^{G}(\phi))(x)=P_{t+s}^{G}(\phi)(x),\  \ s,t\in \lbrack
0,\infty),\ x\in \mathbb{R}.\  \label{Chapman}%
\end{equation}
Such type of nonlinear semigroup was studied in Nisio 1976 \cite{Nisio1},
\cite{Nisio2}.

\begin{proposition}
For each $t>0$, the $G$--normal distribution $P_{t}^{G}$ is a nonlinear
expectation on $\mathcal{H}=lip(\mathbb{R})$, with $\Omega=\mathbb{R}$,
satisfying (a)--(e) of Definition \ref{Def-1}. The corresponding \ completion
space $[\mathcal{H]=[}lip(\mathbb{R})]_{t}$ under the norm $\left \Vert
\phi \right \Vert _{t}:=P_{t}^{G}(|\phi|)(0)$ contains $\phi(x)=x^{n}$,
$n=1,2,\cdots$, as well as $x^{n}\psi$, $\psi \in lip(\mathbb{R)}$ as its
special elements. Relation (\ref{P_tx}) still holds. We also have the
following properties:\newline(1) Central symmetric:\textbf{ }$P_{t}^{G}%
(\phi(\cdot))=P_{t}^{G}(\phi(-\cdot))$;\newline \textbf{(2)} For each convex
$\phi \in \lbrack lip(\mathbb{R})]$ we have
\[
P_{t}^{G}(\phi)(0)=\frac{1}{\sqrt{2\pi t}}\int_{-\infty}^{\infty}\phi
(x)\exp(-\frac{x^{2}}{2t})dx;
\]
For each concave $\phi$, we have, for $\sigma_{0}>0$,%
\[
P_{t}^{G}(\phi)(0)=\frac{1}{\sqrt{2\pi t}\sigma_{0}}\int_{-\infty}^{\infty
}\phi(x)\exp(-\frac{x^{2}}{2t\sigma_{0}^{2}})dx,
\]
and $P_{t}^{G}(\phi)(0)=\phi(0)$ for $\sigma_{0}=0$. In particular, we have
\begin{eqnarray*}
P_{t}^{G}((x)_{x\in \mathbb{R}})  &  =0,\  \  \ P_{t}^{G}((x^{2n+1}%
)_{x\in \mathbb{R}})=P_{t}^{G}((-x^{2n+1})_{x\in \mathbb{R}}),\ n=1,2,\cdots,\\
P_{t}^{G}((x^{2})_{x\in \mathbb{R}})  &  =t,\  \  \ P_{t}^{G}((-x^{2}%
)_{x\in \mathbb{R}})=-\sigma_{0}^{2}t.
\end{eqnarray*}
\end{proposition}

\begin{remark}
Corresponding to the above four expressions, a random $X$ with the $G$--normal
distribution $P_{t}^{G}$ satisfies%
\begin{eqnarray*}
\mathbb{E}[X]  &  =0,\  \  \mathbb{E}[X^{2n+1}]=\mathbb{E}[-X^{2n+1}],\\
\mathbb{E}[X^{2}]  &  =t,\  \  \  \mathbb{E}[-X^{2}]=-\sigma_{0}^{2}t.
\end{eqnarray*}
See the next section for a detail study.
\end{remark}

\section{$1$--dimensional $G$--Brownian motion under $G$--expectation}

In the rest of this paper, we denote by $\Omega=C_{0}(\mathbb{R}^{+})$ the
space of all $\mathbb{R}$--valued continuous paths $(\omega_{t})_{t\in
\mathbb{R}^{+}}$ with $\omega_{0}=0$, equipped with the distance
\[
\rho(\omega^{1},\omega^{2}):=\sum_{i=1}^{\infty}2^{-i}[(\max_{t\in \lbrack
0,i]}|\omega_{t}^{1}-\omega_{t}^{2}|)\wedge1].
\]
We set, for each $t\in \lbrack0,\infty)$,
\begin{align*}
\mathbf{W}_{t}  &  :=\{ \omega_{\cdot \wedge t}:\omega \in \mathbf{\Omega}\},\;
\; \\
\mathcal{F}_{t}  &  :=\mathcal{B}_{t}(\mathbf{W})=\mathcal{B}(\mathbf{W}%
_{t}),\; \; \\
\mathcal{F}_{t+}  &  :=\mathcal{B}_{t+}(\mathbf{W})=\bigcap_{s>t}%
\mathcal{B}_{s}(\mathbf{W}),\\
\mathcal{F}  &  :=%
{\displaystyle \bigvee \limits_{s>t}}
\mathcal{F}_{s}.
\end{align*}
$(\mathbf{\Omega},\mathcal{F})$ is the canonical space equipped with the
natural filtration and $\omega=(\omega_{t})_{t\geq0}$ is the corresponding
canonical process.

For each fixed $T\geq0$, we consider the following space of random variables:
\[
L_{ip}^{0}(\mathcal{F}_{T}):=\{X(\omega)=\phi(\omega_{t_{1}},\cdots
,\omega_{t_{m}}),\forall m\geq1,\;t_{1},\cdots,t_{m}\in \lbrack0,T],\forall
\phi \in lip(\mathbb{R}^{m})\}.
\]
It is clear that $L_{ip}^{0}(\mathcal{F}_{t})\subseteq L_{ip}^{0}%
(\mathcal{F}_{T})$, for $t\leq T$. We also denote%
\[
L_{ip}^{0}(\mathcal{F}):=%
{\displaystyle \bigcup \limits_{n=1}^{\infty}}
L_{ip}^{0}(\mathcal{F}_{n}).
\]

\begin{remark}
It is clear that $lip(\mathbb{R}^{m})$ and then $L_{ip}^{0}(\mathcal{F}_{T})$
and $L_{ip}^{0}(\mathcal{F})$ are vector lattices. Moreover, since $\phi
,\psi \in lip(\mathbb{R}^{m})$ implies $\phi \cdot \psi \in lip(\mathbb{R}^{m})$
thus $X$, $Y\in L_{ip}^{0}(\mathcal{F}_{T})$ implies $X\cdot Y\in L_{ip}%
^{0}(\mathcal{F}_{T})$.
\end{remark}

We will consider the canonical space and set $B_{t}(\omega)=\omega_{t}$,
$t\in \lbrack0,\infty)$, for $\omega \in \Omega$.

\begin{definition}
\label{Def-3}The canonical process $B$ is called a $G$\textbf{--Brownian}
motion under a nonlinear expectation $\mathbb{E}$ defined on $L_{ip}%
^{0}(\mathcal{F})$ if for each $T>0$, $m=1,2,\cdots,$ and for each $\phi \in
lip(\mathbb{R}^{m})$, $0\leq t_{1}<\cdots<t_{m}\leq T$, we have%
\[
\mathbb{E}[\phi(B_{t_{1}},B_{t_{2}}-B_{t_{1}},\cdots,B_{t_{m}}-B_{t_{m-1}%
})]=\phi_{m},
\]
where $\phi_{m}\in \mathbb{R}$ is obtained via the following procedure:%
\begin{eqnarray*}
\phi_{1}(x_{1},\cdots,x_{m-1})  &  =P_{t_{m}-t_{m-1}}^{G}(\phi(x_{1}%
,\cdots,x_{m-1},\cdot));\\
\phi_{2}(x_{1},\cdots,x_{m-2})  &  =P_{t_{m-1}-t_{m-2}}^{G}(\phi_{1}%
(x_{1},\cdots,x_{m-2},\cdot));\\
&  \vdots \\
\phi_{m-1}(x_{1})  &  =P_{t_{2}-t_{1}}^{G}(\phi_{m-2}(x_{1},\cdot));\\
\phi_{m}  &  =P_{t_{1}}^{G}(\phi_{m-1}(\cdot)).
\end{eqnarray*}
The related conditional expectation of $X=\phi(B_{t_{1}},B_{t_{2}}-B_{t_{1}%
},\cdots,B_{t_{m}}-B_{t_{m-1}})$ under $\mathcal{F}_{t_{j}}$ is defined by%
\begin{eqnarray}
\mathbb{E}[X|\mathcal{F}_{t_{j}}]  &  =\mathbb{E}[\phi(B_{t_{1}},B_{t_{2}%
}-B_{t_{1}},\cdots,B_{t_{m}}-B_{t_{m-1}})|\mathcal{F}_{t_{j}}%
]\label{Condition}\\
&  =\phi_{m-j}(B_{t_{1}},\cdots,B_{t_{j}}-B_{t_{j-1}}).\nonumber
\end{eqnarray}
\end{definition}

It is proved in \cite{Peng2005} that $\mathbb{E}[\cdot]$ consistently defines
a nonlinear expectation on the vector lattice $L_{ip}^{0}(\mathcal{F}_{T})$ as
well as on $L_{ip}^{0}(\mathcal{F})$ satisfying (a)--(e) in Definition
\ref{Def-1}. It follows that $\mathbb{E}[|X|]$, $X\in L_{ip}^{0}%
(\mathcal{F}_{T})$ (resp. $L_{ip}^{0}(\mathcal{F})$) forms a norm and that
$L_{ip}^{0}(\mathcal{F}_{T})$ (resp. $L_{ip}^{0}(\mathcal{F})$) can be
continuously extended to a Banach space, denoted by $L_{G}^{1}(\mathcal{F}%
_{T})$ (resp. $L_{G}^{1}(\mathcal{F})$). For each $0\leq t\leq T<\infty$, we
have $L_{G}^{1}(\mathcal{F}_{t})\subseteq L_{G}^{1}(\mathcal{F}_{T})\subset
L_{G}^{1}(\mathcal{F})$. It is easy to check that, in $L_{G}^{1}%
(\mathcal{F}_{T})$ (resp. $L_{G}^{1}(\mathcal{F})$), $\mathbb{E}[\cdot]$ still
satisfies (a)--(e) in Definition \ref{Def-1}.

\begin{definition}
The expectation $\mathbb{E}[\cdot]:L_{G}^{1}(\mathcal{F})\mapsto \mathbb{R}$
introduced through above procedure is called $G$\textbf{--expectation}. The
corresponding canonical process $B$ is called a $G$--Brownian motion under
$\mathbb{E}[\cdot]$.
\end{definition}

For a given $p>1$, we also denote $L_{G}^{p}(\mathcal{F})=\{X\in L_{G}%
^{1}(\mathcal{F}),\ |X|^{p}\in L_{G}^{1}(\mathcal{F})\}$. $L_{G}%
^{p}(\mathcal{F})$ is also a Banach space under the norm $\left \Vert
X\right \Vert _{p}:=(\mathbb{E}[|X|^{p}])^{1/p}$. We have (see Appendix)
\[
\left \Vert X+Y\right \Vert _{p}\leq \left \Vert X\right \Vert _{p}+\left \Vert
Y\right \Vert _{p}%
\]
and, for each $X\in L_{G}^{p}$, $Y\in L_{G}^{q}(Q)$ with $\frac{1}{p}+\frac
{1}{q}=1$,%
\[
\left \Vert XY\right \Vert =\mathbb{E}[|XY|]\leq \left \Vert X\right \Vert
_{p}\left \Vert X\right \Vert _{q}.
\]
With this we have $\left \Vert X\right \Vert _{p}\leq \left \Vert X\right \Vert
_{p^{\prime}}$ if $p\leq p^{\prime}$.

We now consider the conditional expectation introduced in (\ref{Condition}).
For each fixed $t=t_{j}\leq T$, the conditional expectation $\mathbb{E}%
[\cdot|\mathcal{F}_{t}]:L_{ip}^{0}(\mathcal{F}_{T})\mapsto L_{ip}%
^{0}(\mathcal{F}_{t})$ is a continuous mapping under $\left \Vert
\cdot \right \Vert $ since $\mathbb{E}[\mathbb{E}[X|\mathcal{F}_{t}%
]]=\mathbb{E}[X]$, $X\in L_{ip}^{0}(\mathcal{F}_{T})$ and%
\begin{align*}
\mathbb{E}[\mathbb{E}[X|\mathcal{F}_{t}]-\mathbb{E}[Y|\mathcal{F}_{t}]]  &
\leq \mathbb{E}[X-Y],\\
\left \Vert \mathbb{E}[X|\mathcal{F}_{t}]-\mathbb{E}[Y|\mathcal{F}%
_{t}]\right \Vert  &  \leq \left \Vert X-Y\right \Vert .
\end{align*}
It follows that $\mathbb{E}[\cdot|\mathcal{F}_{t}]$ can be also extended as a
continuous mapping $L_{G}^{1}(\mathcal{F}_{T})\mapsto L_{G}^{1}(\mathcal{F}%
_{t})$. If the above $T$ is not fixed, then we can obtain $\mathbb{E}%
[\cdot|\mathcal{F}_{t}]:L_{G}^{1}(\mathcal{F})\mapsto L_{G}^{1}(\mathcal{F}%
_{t})$.

\begin{proposition}
\label{Prop-1-7}We list the properties of $\mathbb{E}[\cdot|\mathcal{F}_{t}]$
that hold in $L_{ip}^{0}(\mathcal{F}_{T})$ and still hold for $X$, $Y\in$
$L_{G}^{1}(\mathcal{F})$:\newline \newline \textbf{(i) }$\mathbb{E}%
[X|\mathcal{F}_{t}]=X$, for $X\in L_{G}^{1}(\mathcal{F}_{t})$, $t\leq
T$.\newline \textbf{(ii) }If $X\geq Y$, then $\mathbb{E}[X|\mathcal{F}_{t}%
]\geq \mathbb{E}[Y|\mathcal{F}_{t}]$.\newline \textbf{(iii) }$\mathbb{E}%
[X|\mathcal{F}_{t}]-\mathbb{E}[Y|\mathcal{F}_{t}]\leq \mathbb{E}%
[X-Y|\mathcal{F}_{t}].$\newline \textbf{(iv) }$\mathbb{E}[\mathbb{E}%
[X|\mathcal{F}_{t}]|\mathcal{F}_{s}]=\mathbb{E}[X|\mathcal{F}_{t\wedge s}]$,
$\mathbb{E}[\mathbb{E}[X|\mathcal{F}_{t}]]=\mathbb{E}[X].$\newline \textbf{(v)}
$\mathbb{E}[X+\eta|\mathcal{F}_{t}]=\mathbb{E}[X|\mathcal{F}_{t}]+\eta$,
$\eta \in L_{G}^{1}(\mathcal{F}_{t}).$\newline \textbf{(vi)} $\mathbb{E}[\eta
X|\mathcal{F}_{t}]=\eta^{+}\mathbb{E}[X|\mathcal{F}_{t}]+\eta^{-}%
\mathbb{E}[-X|\mathcal{F}_{t}]$, for each bounded $\eta \in L_{G}%
^{1}(\mathcal{F}_{t}).$\newline \textbf{(vii)} For each $X\in L_{G}%
^{1}(\mathcal{F}_{T}^{t})$, $\mathbb{E}[X|\mathcal{F}_{t}]=\mathbb{E}%
[X]$,\newline \newline where $L_{G}^{1}(\mathcal{F}_{T}^{t})$ is the extension,
under $\left \Vert \cdot \right \Vert $, of $L_{ip}^{0}(\mathcal{F}_{T}^{t})$
which consists of random variables of the form $\phi(B_{t_{1}}-B_{t_{1}%
},B_{t_{2}}-B_{t_{1}},\cdots,B_{t_{m}}-B_{t_{m-1}})$, $m=1,2,\cdots$, $\phi \in
lip(\mathbb{R}^{m})$, $t_{1},\cdots,t_{m}\in \lbrack t,T]$. Condition (vi) is
the positive homogeneity, see Remark \ref{Rem-Def1}.
\end{proposition}

\begin{definition}
An $X\in L_{G}^{1}(\mathcal{F})$ is said to be independent of $\mathcal{F}%
_{t}$ under the $G$--expectation $\mathbb{E}$ for some given $t\in
\lbrack0,\infty)$, if for each real function $\Phi$ suitably defined on
$\mathbb{R}$ such that $\Phi(X)\in L_{G}^{1}(\mathcal{F})$ we have
\[
\mathbb{E}[\Phi(X)|\mathcal{F}_{t}]=\mathbb{E}[\Phi(X)].
\]
\end{definition}

\begin{remark}
It is clear that all elements in $L_{G}^{1}(\mathcal{F})$ are independent of
$\mathcal{F}_{0}$. Just like the classical situation, the increments of
$G$-Brownian motion $(B_{t+s}-B_{s})_{t\geq0}$ is independent of
$\mathcal{F}_{s}$. In fact it is a new $G$--Brownian motion since, just like
the classical situation, the increments of $B$ are identically distributed.
\end{remark}

\begin{example}
\label{Exam-1}For each $n=0,1,2,\cdots$, $0\leq s-t$, we have $\mathbb{E}%
[B_{t}-B_{s}|\mathcal{F}_{s}]=0$ and, for $n=1,2,\cdots,$
\[
\mathbb{E}[|B_{t}-B_{s}|^{n}|\mathcal{F}_{s}]=\mathbb{E}[|B_{t-s}|^{2n}%
]=\frac{1}{\sqrt{2\pi(t-s)}}\int_{-\infty}^{\infty}|x|^{n}\exp(-\frac{x^{2}%
}{2(t-s)})dx.
\]
But we have%
\[
\mathbb{E}[-|B_{t}-B_{s}|^{n}|\mathcal{F}_{s}]=\mathbb{E}[-|B_{t-s}%
|^{n}]=-\sigma_{0}^{n}\mathbb{E}[|B_{t-s}|^{n}].
\]
Exactly as in classical cases, we have
\begin{align*}
\mathbb{E}[(B_{t}-B_{s})^{2}|\mathcal{F}_{s}]  &  =t-s,\  \  \  \mathbb{E}%
[(B_{t}-B_{s})^{4}|\mathcal{F}_{s}]=3(t-s)^{2},\\
\mathbb{E}[(B_{t}-B_{s})^{6}|\mathcal{F}_{s}]  &  =15(t-s)^{3},\  \  \mathbb{E}%
[(B_{t}-B_{s})^{8}|\mathcal{F}_{s}]=105(t-s)^{4},\\
\mathbb{E}[|B_{t}-B_{s}||\mathcal{F}_{s}]  &  =\frac{\sqrt{2(t-s)}}{\sqrt{\pi
}},\  \  \mathbb{E}[|B_{t}-B_{s}|^{3}|\mathcal{F}_{s}]=\frac{2\sqrt
{2}(t-s)^{3/2}}{\sqrt{\pi}},\\
\mathbb{E}[|B_{t}-B_{s}|^{5}|\mathcal{F}_{s}]  &  =8\frac{\sqrt{2}(t-s)^{5/2}%
}{\sqrt{\pi}}.
\end{align*}
\end{example}

\begin{example}
\label{Exam-2}For each $n=1,2,\cdots,$ $0\leq s\leq t<T$ and $X\in L_{G}%
^{1}(\mathcal{F}_{s})$, since $\mathbb{E[}B_{T-t}^{2n-1}]=\mathbb{E[-}%
B_{T-t}^{2n-1}]$, we have, by (vi) of Proposition \ref{Prop-1-7},
\begin{align*}
\mathbb{E}[X(B_{T}-B_{t})^{2n-1}]  &  =\mathbb{E}[X^{+}\mathbb{E[}(B_{T}%
-B_{t})^{2n-1}|\mathcal{F}_{t}]+X^{-}\mathbb{E[-}(B_{T}-B_{t})^{2n-1}%
|\mathcal{F}_{t}]]\\
&  =\mathbb{E}[|X|]\cdot \mathbb{E[}B_{T-t}^{2n-1}],\\
\mathbb{E}[X(B_{T}-B_{t})|\mathcal{F}_{s}]  &  =\mathbb{E}[-X(B_{T}%
-B_{t})|\mathcal{F}_{s}]=0.
\end{align*}
We also have%
\[
\mathbb{E}[X(B_{T}-B_{t})^{2}|\mathcal{F}_{t}]=X^{+}(T-t)-\sigma_{0}^{2}%
X^{-}(T-t).
\]
\end{example}

\begin{remark}
It is clear that we can define an expectation $E[\cdot]$ on $L_{ip}%
^{0}(\mathcal{F})$ in the same way as in Definition \ref{Def-3} with the
standard normal distribution $P_{1}(\cdot)$ in the place of $P_{1}^{G}(\cdot
)$. Since $P_{1}(\cdot)$ is dominated by $P_{1}^{G}(\cdot)$ in the sense
$P_{1}(\phi)-P_{1}(\psi)\leq P_{1}^{G}(\phi-\psi)$, then $E[\cdot]$ can be
continuously extended to $L_{G}^{1}(\mathcal{F})$. $E[\cdot]$ is a linear
expectation under which $(B_{t})_{t\geq0}$ behaves as a Brownian motion. We
have
\begin{equation}
E[X]\leq \mathbb{E}[X],\  \  \forall X\in L_{G}^{1}(\mathcal{F}). \label{Eg-domi}%
\end{equation}
In particular, $\mathbb{E[}B_{T-t}^{2n-1}]=\mathbb{E[-}B_{T-t}^{2n-1}]\geq
E\mathbb{[-}B_{T-t}^{2n-1}]=0$. Such kind of extension under a domination
relation was discussed in details in \cite{Peng2005}.
\end{remark}

The following property is very useful

\begin{proposition}
\label{E-x+y}Let $X,Y\in L_{G}^{1}(\mathcal{F})$ be such that $\mathbb{E}%
[Y]=-\mathbb{E}[-Y]$ (thus $\mathbb{E}[Y]=E[Y]$), then we have%
\[
\mathbb{E}[X+Y]=\mathbb{E}[X]+\mathbb{E}[Y].
\]
In particular, if $\mathbb{E}[Y]=\mathbb{E}[-Y]=0$, then $\mathbb{E}%
[X+Y]=\mathbb{E}[X]$.
\end{proposition}

\begin{proof}
It is simply because we have $\mathbb{E}[X+Y]\leq \mathbb{E}[X]+\mathbb{E}[Y]$
and
\[
\mathbb{E}[X+Y]\geq \mathbb{E}[X]-\mathbb{E}[-Y]=\mathbb{E}[X]+\mathbb{E}%
[Y]\text{.}%
\]
\end{proof}

\begin{example}
\label{Exam-B2}We have%
\begin{align*}
\mathbb{E}[B_{t}^{2}-B_{s}^{2}|\mathcal{F}_{s}]  &  =\mathbb{E}[(B_{t}%
-B_{s}+B_{s})^{2}-B_{s}^{2}|\mathcal{F}_{s}]\\
&  =E[(B_{t}-B_{s})^{2}+2(B_{t}-B_{s})B_{s}|\mathcal{F}_{s}]\\
&  =t-s,
\end{align*}
since $2(B_{t}-B_{s})B_{s}$ satisfies the condition for $Y$ in Proposition
\ref{E-x+y}, and%
\begin{align*}
\mathbb{E}[(B_{t}^{2}-B_{s}^{2})^{2}|\mathcal{F}_{s}]  &  =\mathbb{E}%
[\{(B_{t}-B_{s}+B_{s})^{2}-B_{s}^{2}\}^{2}|\mathcal{F}_{s}]\\
&  =\mathbb{E}[\{(B_{t}-B_{s})^{2}+2(B_{t}-B_{s})B_{s}\}^{2}|\mathcal{F}%
_{s}]\\
&  =\mathbb{E}[(B_{t}-B_{s})^{4}+4(B_{t}-B_{s})^{3}B_{s}+4(B_{t}-B_{s}%
)^{2}B_{s}^{2}|\mathcal{F}_{s}]\\
&  \leq \mathbb{E}[(B_{t}-B_{s})^{4}]+4\mathbb{E}[|B_{t}-B_{s}|^{3}%
]|B_{s}|+4(t-s)B_{s}^{2}\\
&  =3(t-s)^{2}+8(t-s)^{3/2}|B_{s}|+4(t-s)B_{s}^{2}.
\end{align*}
\end{example}

\section{It\^{o}'s integral of $G$--Brownian motion}

\subsection{Bochner's integral}

\begin{definition}
\label{Def-4}For $T\in \mathbb{R}_{+}$, a partition $\pi_{T}$ of $[0,T]$ is a
finite ordered subset $\pi=\{t_{1},\cdots,t_{N}\}$ such that $0=t_{0}%
<t_{1}<\cdots<t_{N}=T$. We denote
\[
\mu(\pi_{T})=\max \{|t_{i+1}-t_{i}|,i=0,1,\cdots,N-1\} \text{.}%
\]
We use $\pi_{T}^{N}=\{t_{0}^{N}<t_{1}^{N}<\cdots<t_{N}^{N}\}$ to denote a
sequence of partitions of $[0,T]$ such that $\lim_{N\rightarrow \infty}\mu
(\pi_{T}^{N})=0$.
\end{definition}

Let $p\geq1$ be fixed. We consider the following type of simple processes: for
a given partition $\{t_{0},\cdots,t_{N}\}=\pi_{T}$ of $[0,T]$, we set%
\[
\eta_{t}(\omega)=\sum_{j=0}^{N-1}\xi_{j}(\omega)\mathbf{I}_{[t_{j},t_{j+1}%
)}(t),
\]
where $\xi_{i}\in L_{G}^{p}(\mathcal{F}_{t_{i}})$, $i=0,1,2,\cdots,N-1$, are
given. The collection and these type of processes is denoted by $M_{G}%
^{p,0}(0,T)$.

\begin{definition}
\label{Def-5}For an $\eta \in M_{G}^{1,0}(0,T)$ with $\eta_{t}=\sum_{j=0}%
^{N-1}\xi_{j}(\omega)\mathbf{I}_{[t_{j},t_{j+1})}(t)$, the related Bochner
integral is
\[
\int_{0}^{T}\eta_{t}(\omega)dt=\sum_{j=0}^{N-1}\xi_{j}(\omega)(t_{j+1}%
-t_{j}).
\]
\end{definition}

\begin{remark}
We set, for each $\eta \in M_{G}^{1,0}(0,T)$,
\[
\mathbb{\tilde{E}}_{T}[\eta]:=\frac{1}{T}\int_{0}^{T}\mathbb{E}[\eta
_{t}]dt=\frac{1}{T}\sum_{j=0}^{N-1}\mathbb{E[}\xi_{j}(\omega)](t_{j+1}%
-t_{j}).
\]
It is easy to check that $\mathbb{\tilde{E}}_{T}:M_{G}^{1,0}(0,T)\longmapsto
\mathbb{R}$ forms a nonlinear expectation satisfying (a)--(e) of Definition
\ref{Def-1}. By Remark \ref{Rem-2}, we can introduce a natural norm
$\left \Vert \eta \right \Vert _{T}^{1}=\mathbb{\tilde{E}}_{T}[|\eta|]=\frac
{1}{T}\int_{0}^{T}\mathbb{E}[|\eta_{t}|]dt$. Under this norm $M_{G}%
^{1,0}(0,T)$ can be continuously extended to $M_{G}^{1}(0,T)$ which is a
Banach space.
\end{remark}

\begin{definition}
For each $p\geq1$, we will denote by $M_{G}^{p}(0,T)$ the completion of
$M_{G}^{p,0}(0,T)$ under the norm%
\[
(\frac{1}{T}\int_{0}^{T}\left \Vert \eta_{t}^{p}\right \Vert dt)^{1/p}=\left(
\frac{1}{T}\sum_{j=0}^{N-1}\mathbb{E[}|\xi_{j}(\omega)|^{p}](t_{j+1}%
-t_{j})\right)  ^{1/p}.
\]
\end{definition}

We observe that,
\[
\mathbb{E}[|\int_{0}^{T}\eta_{t}(\omega)dt|]\leq \sum_{j=0}^{N-1}\left \Vert
\xi_{j}(\omega)\right \Vert (t_{j+1}-t_{j})=\int_{0}^{T}\mathbb{E}[|\eta
_{t}|]dt.
\]
We then have

\begin{proposition}
The linear mapping $\int_{0}^{T}\eta_{t}(\omega)dt:M_{G}^{1,0}(0,T)\mapsto
L_{G}^{1}(\mathcal{F}_{T})$ is continuous. and thus can be continuously
extended to $M_{G}^{1}(0,T)\mapsto L_{G}^{1}(\mathcal{F}_{T})$. We still
denote this extended mapping by $\int_{0}^{T}\eta_{t}(\omega)dt$, $\eta \in
M_{G}^{1}(0,T)$. We have%
\begin{equation}
\mathbb{E}[|\int_{0}^{T}\eta_{t}(\omega)dt|]\leq \int_{0}^{T}\mathbb{E}%
[|\eta_{t}|]dt,\  \  \  \forall \eta \in M_{G}^{1}(0,T). \label{ine-dt}%
\end{equation}
\end{proposition}

Since $M_{G}^{1}(0,T)\supset M_{G}^{p}(0,T)$, for $p\geq1$, this definition
holds for $\eta \in M_{G}^{p}(0,T)$.

\subsection{It\^{o}'s integral of $G$--Brownian motion}

\begin{definition}
For each $\eta \in M_{G}^{2,0}(0,T)$ with the form $\eta_{t}(\omega)=\sum
_{j=0}^{N-1}\xi_{j}(\omega)\mathbf{I}_{[t_{j},t_{j+1})}(t)$, we define
\[
I(\eta)=\int_{0}^{T}\eta(s)dB_{s}:=\sum_{j=0}^{N-1}\xi_{j}(B_{t_{j+1}%
}-B_{t_{j}})\mathbf{.}%
\]
\end{definition}

\begin{lemma}
\label{bdd}The mapping $I:M_{G}^{2,0}(0,T)\longmapsto L_{G}^{2}(\mathcal{F}%
_{T})$ is a linear continuous mapping and thus can be continuously extended to
$I:M_{G}^{2}(0,T)\longmapsto L_{G}^{2}(\mathcal{F}_{T})$. In fact we have
\begin{align}
\mathbb{E}[\int_{0}^{T}\eta(s)dB_{s}]  &  =0,\  \  \label{e1}\\
\mathbb{E}[(\int_{0}^{T}\eta(s)dB_{s})^{2}]  &  \leq \int_{0}^{T}%
\mathbb{E}[(\eta(t))^{2}]dt. \label{e2}%
\end{align}
\end{lemma}

\begin{definition}
We define, for a fixed $\eta \in M_{G}^{2}(0,T)$, the stochastic integral
\[
\int_{0}^{T}\eta(s)dB_{s}:=I(\eta).
\]
It is clear that (\ref{e1}), (\ref{e2}) still hold for $\eta \in M_{G}%
^{2}(0,T)$.
\end{definition}

\textbf{Proof of Lemma \ref{bdd}. }From Example \ref{Exam-2}, for each $j$,
\[
\mathbb{E}\mathbf{[}\xi_{j}(B_{t_{j+1}}-B_{t_{j}})|\mathcal{F}_{t_{j}}]=0.
\]
We have%
\begin{align*}
\mathbb{E}[\int_{0}^{T}\eta(s)dB_{s}]  &  =\mathbb{E[}\int_{0}^{t_{N-1}}%
\eta(s)dB_{s}+\xi_{N-1}(B_{t_{N}}-B_{t_{N-1}})]\\
&  =\mathbb{E[}\int_{0}^{t_{N-1}}\eta(s)dB_{s}+\mathbb{E}\mathbf{[}\xi
_{N-1}(B_{t_{N}}-B_{t_{N-1}})|\mathcal{F}_{t_{N-1}}]]\\
&  =\mathbb{E[}\int_{0}^{t_{N-1}}\eta(s)dB_{s}].
\end{align*}
We then can repeat this procedure to obtain (\ref{e1}). We now prove
(\ref{e2}):
\begin{align*}
\mathbb{E}[(\int_{0}^{T}\eta(s)dB_{s})^{2}]  &  =\mathbb{E[}\left(  \int
_{0}^{t_{N-1}}\eta(s)dB_{s}+\xi_{N-1}(B_{t_{N}}-B_{t_{N-1}})\right)  ^{2}]\\
&  =\mathbb{E[}\left(  \int_{0}^{t_{N-1}}\eta(s)dB_{s}\right)  ^{2}\\
&  +\mathbb{E}[2\left(  \int_{0}^{t_{N-1}}\eta(s)dB_{s}\right)  \xi
_{N-1}(B_{t_{N}}-B_{t_{N-1}})+\xi_{N-1}^{2}(B_{t_{N}}-B_{t_{N-1}}%
)^{2}|\mathcal{F}_{t_{N-1}}]]\\
&  =\mathbb{E[}\left(  \int_{0}^{t_{N-1}}\eta(s)dB_{s}\right)  ^{2}+\xi
_{N-1}^{2}(t_{N}-t_{N-1})].
\end{align*}
Thus $\mathbb{E}[(\int_{0}^{t_{N}}\eta(s)dB_{s})^{2}]\leq \mathbb{E[}\left(
\int_{0}^{t_{N-1}}\eta(s)dB_{s}\right)  ^{2}]+\mathbb{E}[\xi_{N-1}^{2}%
](t_{N}-t_{N-1})]$. We then repeat this procedure to deduce
\[
\mathbb{E}[(\int_{0}^{T}\eta(s)dB_{s})^{2}]\leq \sum_{j=0}^{N-1}\mathbb{E}%
[(\xi_{j})^{2}](t_{j+1}-t_{j})=\int_{0}^{T}\mathbb{E}[(\eta(t))^{2}]dt.
\]
$\blacksquare$

We list some main properties of the It\^{o}'s integral of $G$--Brownian
motion. We denote for some $0\leq s\leq t\leq T$,
\[
\int_{s}^{t}\eta_{u}dB_{u}:=\int_{0}^{T}\mathbf{I}_{[s,t]}(u)\eta_{u}dB_{u}.
\]
We have

\begin{proposition}
\label{Prop-Integ}Let $\eta,\theta \in M_{G}^{2}(0,T)$ and let $0\leq s\leq
r\leq t\leq T$. Then in $L_{G}^{1}(\mathcal{F}_{T})$ we have\newline(i)
$\int_{s}^{t}\eta_{u}dB_{u}=\int_{s}^{r}\eta_{u}dB_{u}+\int_{r}^{t}\eta
_{u}dB_{u},$\newline(ii) $\int_{s}^{t}(\alpha \eta_{u}+\theta_{u})dB_{u}%
=\alpha \int_{s}^{t}\eta_{u}dB_{u}+\int_{s}^{t}\theta_{u}dB_{u},\ $if$\  \alpha$
is bounded and in $L_{G}^{1}(\mathcal{F}_{s})$,\newline(iii) $\mathbb{E[}%
X+\int_{r}^{T}\eta_{u}dB_{u}|\mathcal{F}_{s}]=\mathbb{E[}X]$, $\forall X\in
L_{G}^{1}(\mathcal{F})$.
\end{proposition}

\subsection{Quadratic variation process of $G$--Brownian motion}

We now study a very interesting process of the $G$-Brownian motion. Let
$\pi_{t}^{N}$, $N=1,2,\cdots$, be a sequence of partitions of $[0,t]$. We consider%

\begin{align*}
B_{t}^{2}  &  =\sum_{j=0}^{N-1}[B_{t_{j+1}^{N}}^{2}-B_{t_{j}^{N}}^{2}]\\
&  =\sum_{j=0}^{N-1}2B_{t_{j}^{N}}(B_{t_{j+1}^{N}}-B_{t_{j}^{N}})+\sum
_{j=0}^{N-1}(B_{t_{j+1}^{N}}-B_{t_{j}^{N}})^{2}.
\end{align*}
As $\mu(\pi_{t}^{N})\rightarrow0$, the first term of the right side tends to
$\int_{0}^{t}B_{s}dB_{s}$. The second term must converge. We denote its limit
by $\left \langle B\right \rangle _{t}$, i.e.,
\begin{equation}
\left \langle B\right \rangle _{t}=\lim_{\mu(\pi_{t}^{N})\rightarrow0}\sum
_{j=0}^{N-1}(B_{t_{j+1}^{N}}-B_{t_{j}^{N}})^{2}=B_{t}^{2}-2\int_{0}^{t}%
B_{s}dB_{s}. \label{quadra-def}%
\end{equation}
By the above construction, $\left \langle B\right \rangle _{t}$, $t\geq0$, is an
increasing process with $\left \langle B\right \rangle _{0}=0$. We call it the
\textbf{quadratic variation process} of the $G$--Brownian motion $B$. Clearly
$\left \langle B\right \rangle $ is an increasing process. It perfectly
characterizes the part of uncertainty, or ambiguity, of $G$--Brownian motion.
It is important to keep in mind that $\left \langle B\right \rangle _{t}$ is not
a deterministic process unless the case $\sigma=1$, i.e., when $B$ is a
classical Brownian motion. In fact we have

\begin{lemma}
\label{Lem-Q1}We have, for each $0\leq s\leq t<\infty$%
\begin{align}
\mathbb{E}[\left \langle B\right \rangle _{t}-\left \langle B\right \rangle
_{s}|\mathcal{F}_{s}]  &  =t-s,\  \  \label{quadra}\\
\mathbb{E}[-(\left \langle B\right \rangle _{t}-\left \langle B\right \rangle
_{s})|\mathcal{F}_{s}]  &  =-\sigma_{0}^{2}(t-s). \label{quadra1}%
\end{align}
\end{lemma}

\begin{proof}
By the definition of $\left \langle B\right \rangle $ and Proposition
\ref{Prop-Integ}-(iii),
\begin{align*}
\mathbb{E}[\left \langle B\right \rangle _{t}-\left \langle B\right \rangle
_{s}|\mathcal{F}_{s}]  &  =\mathbb{E}[B_{t}^{2}-B_{s}^{2}-2\int_{s}^{t}%
B_{u}dB_{u}|\mathcal{F}_{s}]\\
&  =\mathbb{E}[B_{t}^{2}-B_{s}^{2}|\mathcal{F}_{s}]=t-s.
\end{align*}
The last step can be check as in Example \ref{Exam-B2}. We then have
(\ref{quadra}). (\ref{quadra1}) can be proved analogously with the
consideration of $\mathbb{E}[-(B_{t}^{2}-B_{s}^{2})|\mathcal{F}_{s}%
]=-\sigma^{2}(t-s)$.
\end{proof}

To define the integration of a process $\eta \in M_{G}^{1}(0,T)$ with respect
to $d\left \langle B\right \rangle $, we first define a mapping:%
\[
Q_{0,T}(\eta)=\int_{0}^{T}\eta(s)d\left \langle B\right \rangle _{s}:=\sum
_{j=0}^{N-1}\xi_{j}(\left \langle B\right \rangle _{t_{j+1}}-\left \langle
B\right \rangle _{t_{j}}):M_{G}^{1,0}(0,T)\mapsto L^{1}(\mathcal{F}_{T}).
\]

\begin{lemma}
\label{Lem-Q2}For each $\eta \in M_{G}^{1,0}(0,T)$,
\begin{equation}
\mathbb{E}[|Q_{0,T}(\eta)|]\leq \int_{0}^{T}\mathbb{E}[|\eta_{s}%
|]ds.\  \label{dA}%
\end{equation}
Thus $Q_{0,T}:M_{G}^{1,0}(0,T)\mapsto L^{1}(\mathcal{F}_{T})$ is a continuous
linear mapping. Consequently, $Q_{0,T}$ can be uniquely extended to
$L_{\mathcal{F}}^{1}(0,T)$. We still denote this mapping \ by%
\[
\int_{0}^{T}\eta(s)d\left \langle B\right \rangle _{s}=Q_{0,T}(\eta),\  \  \eta \in
M_{G}^{1}(0,T)\text{.}%
\]
We still have
\begin{equation}
\mathbb{E}[|\int_{0}^{T}\eta(s)d\left \langle B\right \rangle _{s}|]\leq \int
_{0}^{T}\mathbb{E}[|\eta_{s}|]ds,\  \  \forall \eta \in M_{G}^{1}(0,T)\text{.}
\label{qua-ine}%
\end{equation}
\end{lemma}

\begin{proof}
By applying Lemma \ref{Lem-Q1}, (\ref{dA}) can be checked as follows:%
\begin{align*}
\mathbb{E}[|\sum_{j=0}^{N-1}\xi_{j}(\left \langle B\right \rangle _{t_{j+1}%
}-\left \langle B\right \rangle _{t_{j}})|]  &  \leq \sum_{j=0}^{N-1}%
\mathbb{E[}|\xi_{j}|\cdot \mathbb{E}[\left \langle B\right \rangle _{t_{j+1}%
}-\left \langle B\right \rangle _{t_{j}}|\mathcal{F}_{t_{j}}]]\\
&  =\sum_{j=0}^{N-1}\mathbb{E[}|\xi_{j}|](t_{j+1}-t_{j})\\
&  =\int_{0}^{T}\mathbb{E}[|\eta_{s}|]ds.
\end{align*}
\end{proof}

A very interesting point of the quadratic variation process $\left \langle
B\right \rangle $ is, just like the $G$--Brownian motion $B$ it's self, the
increment $\left \langle B\right \rangle _{t+s}-\left \langle B\right \rangle
_{s}$ is independent of $\mathcal{F}_{s}$ and identically distributed like
$\left \langle B\right \rangle _{t}$. In fact we have

\begin{lemma}
\label{Lem-Qua2}For each fixed $s\geq0$, $(\left \langle B\right \rangle
_{s+t}-\left \langle B\right \rangle _{s})_{t\geq0}$ is independent of
$\mathcal{F}_{s}$. It is the quadratic variation process of the Brownian
motion $B_{t}^{s}=B_{s+t}-B_{s}$, $t\geq0$, i.e., $\left \langle B\right \rangle
_{s+t}-\left \langle B\right \rangle _{s}=\left \langle B^{s}\right \rangle _{t}$.
We have
\begin{equation}
\mathbb{E}[\left \langle B^{s}\right \rangle _{t}^{2}|\mathcal{F}_{s}%
]=\mathbb{E}[\left \langle B\right \rangle _{t}^{2}]=t^{2} \label{Qua2}%
\end{equation}
as well as%
\[
\mathbb{E}[\left \langle B^{s}\right \rangle _{t}^{3}|\mathcal{F}_{s}%
]=\mathbb{E}[\left \langle B\right \rangle _{t}^{2}]=t^{3},\  \  \  \mathbb{E}%
[\left \langle B^{s}\right \rangle _{t}^{4}|\mathcal{F}_{s}]=\mathbb{E}%
[\left \langle B\right \rangle _{t}^{4}]=t^{4}.
\]
\end{lemma}

\begin{proof}
The independence is simply from
\begin{align*}
\left \langle B\right \rangle _{s+t}-\left \langle B\right \rangle _{s}  &
=B_{t+s}^{2}-2\int_{0}^{s+t}B_{r}dB_{r}-[B_{s}^{2}-2\int_{0}^{s}B_{r}dB_{r}]\\
&  =(B_{t+s}-B_{s})^{2}-2\int_{s}^{s+t}(B_{r}-B_{s})d(B_{r}-B_{s})\\
&  =\left \langle B^{s}\right \rangle _{t}.
\end{align*}
We set $\phi(t):=\mathbb{E}[\left \langle B\right \rangle _{t}^{2}]$.
\begin{align*}
\phi(t)  &  =\mathbb{E}[\{(B_{t})^{2}-2\int_{0}^{t}B_{u}dB_{u}\}^{2}]\\
&  \leq2\mathbb{E}[(B_{t})^{4}]+8\mathbb{E}[(\int_{0}^{t}B_{u}dB_{u})^{2}]\\
&  \leq6t^{2}+8\int_{0}^{t}\mathbb{E[(}B_{u})^{2}]du\\
&  =10t^{2}.
\end{align*}
This also implies $\mathbb{E}[(\left \langle B\right \rangle _{t+s}-\left \langle
B\right \rangle _{s})^{2}]=\phi(t)\leq14t$. Thus
\begin{align*}
\phi(t)  &  =\mathbb{E}[\{ \left \langle B\right \rangle _{s}+\left \langle
B\right \rangle _{s+t}-\left \langle B\right \rangle _{s}\}^{2}]\\
&  \leq \mathbb{E}[(\left \langle B\right \rangle _{s})^{2}]+\mathbb{E}%
[(\left \langle B^{s}\right \rangle _{t})^{2}]+2\mathbb{E}[\left \langle
B\right \rangle _{s}\left \langle B^{s}\right \rangle _{t}]\\
&  =\phi(s)+\phi(t)+2\mathbb{E}[\left \langle B\right \rangle _{s}%
\mathbb{E}[\left \langle B^{s}\right \rangle _{t}]]\\
&  =\phi(s)+\phi(t)+2st.
\end{align*}
We set $\delta_{N}=t/N$, $t_{k}^{N}=kt/N=k\delta_{N}$ for a positive integer
$N$. By the above inequalities%
\begin{align*}
\phi(t_{N}^{N})  &  \leq \phi(t_{N-1}^{N})+\phi(\delta_{N})+2t_{N-1}^{N}%
\delta_{N}\\
&  \leq \phi(t_{N-2}^{N})+2\phi(\delta_{N})+2(t_{N-1}^{N}+t_{N-2}^{N}%
)\delta_{N}\\
&  \vdots
\end{align*}
We then have
\[
\phi(t)\leq N\phi(\delta_{N})+2\sum_{k=0}^{N-1}t_{k}^{N}\delta_{N}\leq
10\frac{t^{2}}{N}+2\sum_{k=0}^{N-1}t_{k}^{N}\delta_{N}.
\]
Let $N\rightarrow \infty$ we have $\phi(t)\leq2\int_{0}^{t}sds=t^{2}$. Thus
$\mathbb{E}[\left \langle B_{t}\right \rangle ^{2}]\leq t^{2}$. This with
$\mathbb{E}[\left \langle B_{t}\right \rangle ^{2}]\geq E[\left \langle
B_{t}\right \rangle ^{2}]=t^{2}$ implies (\ref{Qua2}).
\end{proof}

\begin{proposition}
\label{Prop-temp}Let $0\leq s\leq t$, $\xi \in L_{G}^{1}(\mathcal{F}_{s})$.
Then%
\begin{align*}
\mathbb{E}[X+\xi(B_{t}^{2}-B_{s}^{2})]  &  =\mathbb{E}[X+\xi(B_{t}-B_{s}%
)^{2}]\\
&  =\mathbb{E}[X+\xi(\left \langle B\right \rangle _{t}-\left \langle
B\right \rangle _{s})].
\end{align*}
\end{proposition}

\begin{proof}
By (\ref{quadra-def}) and Proposition \ref{E-x+y}, we have%
\begin{align*}
\mathbb{E}[X+\xi(B_{t}^{2}-B_{s}^{2})]  &  =\mathbb{E}[X+\xi(\left \langle
B\right \rangle _{t}-\left \langle B\right \rangle _{s}+2\int_{s}^{t}B_{u}%
dB_{u})]\\
&  =\mathbb{E}[X+\xi(\left \langle B\right \rangle _{t}-\left \langle
B\right \rangle _{s})].
\end{align*}
We also have
\begin{align*}
\mathbb{E}[X+\xi(B_{t}^{2}-B_{s}^{2})]  &  =\mathbb{E}[X+\xi \{(B_{t}%
-B_{s})^{2}+2(B_{t}-B_{s})B_{s}^{{}}\}]\\
&  =\mathbb{E}[X+\xi(B_{t}-B_{s})^{2}].
\end{align*}
\end{proof}

We have the following isometry

\begin{proposition}
Let $\eta \in M_{G}^{2}(0,T)$. We have%
\begin{equation}
\mathbb{E}[(\int_{0}^{T}\eta(s)dB_{s})^{2}]=\mathbb{E}[\int_{0}^{T}\eta
^{2}(s)d\left \langle B\right \rangle _{s}]. \label{isometry}%
\end{equation}
\end{proposition}

\begin{proof}
We first consider $\eta \in M_{G}^{2,0}(0,T)$ with the form
\[
\eta_{t}(\omega)=\sum_{j=0}^{N-1}\xi_{j}(\omega)\mathbf{I}_{[t_{j},t_{j+1}%
)}(t)
\]
and thus $\int_{0}^{T}\eta(s)dB_{s}:=\sum_{j=0}^{N-1}\xi_{j}(B_{t_{j+1}%
}-B_{t_{j}})$\textbf{.} By Proposition \ref{E-x+y} we have
\[
\mathbb{E}[X+2\xi_{j}(B_{t_{j+1}}-B_{t_{j}})\xi_{i}(B_{t_{i+1}}-B_{t_{i}%
})]=\mathbb{E}[X]\text{, for }X\in L_{G}^{1}(\mathcal{F)}\text{, }i\not =j.
\]
Thus%
\[
\mathbb{E}[(\int_{0}^{T}\eta(s)dB_{s})^{2}]=\mathbb{E[}\left(  \sum
_{j=0}^{N-1}\xi_{j}(B_{t_{j+1}}-B_{t_{j}})\right)  ^{2}]=\mathbb{E[}\sum
_{j=0}^{N-1}\xi_{j}^{2}(B_{t_{j+1}}-B_{t_{j}})^{2}].
\]
This with Proposition \ref{Prop-temp}, it follows that
\[
\mathbb{E}[(\int_{0}^{T}\eta(s)dB_{s})^{2}]=\mathbb{E[}\sum_{j=0}^{N-1}\xi
_{j}^{2}(\left \langle B\right \rangle _{t_{j+1}}-\left \langle B\right \rangle
_{t_{j}})]=\mathbb{E[}\int_{0}^{T}\eta^{2}(s)d\left \langle B\right \rangle
_{s}].
\]
Thus (\ref{isometry}) holds for $\eta \in M_{G}^{2,0}(0,T)$. We thus can
continuously extend the above equality to the case $\eta \in M_{G}^{2}(0,T)$
and prove (\ref{isometry}).
\end{proof}

\subsection{It\^{o}'s formula for $G$--Brownian motion}

We have the corresponding It\^{o}'s formula of $\Phi(X_{t})$ for a
\textquotedblleft$G$-It\^{o} process\textquotedblright \ $X$. For
simplification, we only treat the case where the function $\Phi$ is
sufficiently regular. We first consider a simple situation.

\label{Lem-26}Let $\Phi \in C^{2}(\mathbb{R}^{n})$ be bounded with bounded
derivatives and $\{ \partial_{x^{\mu}x^{\nu}}^{2}\Phi \}_{\mu,\nu=1}^{n}$ are
uniformly Lipschitz. Let $s\in \lbrack0,T]$ be fixed and let $X=(X^{1}%
,\cdots,X^{n})^{T}$ be an $n$--dimensional process on $[s,T]$ of the form
\[
X_{t}^{\nu}=X_{s}^{\nu}+\alpha^{\nu}(t-s)+\eta^{\nu}(\left \langle
B\right \rangle _{t}-\left \langle B\right \rangle _{s})+\beta^{\nu}(B_{t}%
-B_{s}),
\]
where, for $\nu=1,\cdots,n$, $\alpha^{\nu}$, $\eta^{\nu}$ and $\beta^{\nu}$,
are bounded elements of $L_{G}^{2}(\mathcal{F}_{s})$ and $X_{s}=(X_{s}%
^{1},\cdots,X_{s}^{n})^{T}$ is a given $\mathbb{R}^{n}$--vector in $L_{G}%
^{2}(\mathcal{F}_{s})$. Then we have
\begin{align}
\Phi(X_{t})-\Phi(X_{s})  &  =\int_{s}^{t}\partial_{x^{\nu}}\Phi(X_{u}%
)\beta^{\nu}dB_{u}+\int_{s}^{t}\partial_{x_{\nu}}\Phi(X_{u})\alpha^{\nu
}du\label{B-Ito}\\
&  +\int_{s}^{t}[D_{x^{\nu}}\Phi(X_{u})\eta^{\nu}+\frac{1}{2}\partial_{x^{\mu
}x^{\nu}}^{2}\Phi(X_{u})\beta^{\mu}\beta^{\nu}]d\left \langle B\right \rangle
_{u}.\nonumber
\end{align}
Here we use the Einstein convention, i.e., each single term with repeated
indices $\mu$ and/or $\nu$ implies the summation.

\begin{proof}
For each positive integer $N$ we set $\delta=(t-s)/N$ and take the partition
\[
\pi_{\lbrack s,t]}^{N}=\{t_{0}^{N},t_{1}^{N},\cdots,t_{N}^{N}\}=\{s,s+\delta
,\cdots,s+N\delta=t\}.
\]
We have
\begin{align}
\Phi(X_{t})  &  =\Phi(X_{s})+\sum_{k=0}^{N-1}[\Phi(X_{t_{k+1}^{N}}%
)-\Phi(X_{t_{k}^{N}})]\nonumber \\
&  =\Phi(X_{s})+\sum_{k=0}^{N-1}[\partial_{x^{\mu}}\Phi(X_{t_{k}^{N}%
})(X_{t_{k+1}^{N}}^{\mu}-X_{t_{k}^{N}}^{\mu})\nonumber \\
&  +\frac{1}{2}[\partial_{x^{\mu}x^{\nu}}^{2}\Phi(X_{t_{k}^{N}})(X_{t_{k+1}%
^{N}}^{\mu}-X_{t_{k}^{N}}^{\mu})(X_{t_{k+1}^{N}}^{\nu}-X_{t_{k}^{N}}^{\nu
})+\eta_{k}^{N}]] \label{Ito}%
\end{align}
where
\[
\eta_{k}^{N}=[\partial_{x^{\mu}x^{\nu}}^{2}\Phi(X_{t_{k}^{N}}+\theta
_{k}(X_{t_{k+1}^{N}}-X_{t_{k}^{N}}))-\partial_{x^{\mu}x^{\nu}}^{2}%
\Phi(X_{t_{k}^{N}})](X_{t_{k+1}^{N}}^{\mu}-X_{t_{k}^{N}}^{\mu})(X_{t_{k+1}%
^{N}}^{\nu}-X_{t_{k}^{N}}^{\nu})
\]
with $\theta_{k}\in \lbrack0,1]$. We have%
\begin{align*}
\mathbb{E}[|\eta_{k}^{N}|]  &  =\mathbb{E}[|[\partial_{x^{\mu}x^{\nu}}^{2}%
\Phi(X_{t_{k}^{N}}+\theta_{k}(X_{t_{k+1}^{N}}-X_{t_{k}^{N}}))-\partial
_{x^{\mu}x^{\nu}}^{2}\Phi(X_{t_{k}^{N}})](X_{t_{k+1}^{N}}^{\mu}-X_{t_{k}^{N}%
}^{\mu})(X_{t_{k+1}^{N}}^{\nu}-X_{t_{k}^{N}}^{\nu})|]\\
&  \leq c\mathbb{E[}|X_{t_{k+1}^{N}}-X_{t_{k}^{N}}|^{3}]\leq C[\delta
^{3}+\delta^{3/2}],
\end{align*}
where $c$ is the Lipschitz constant of $\{ \partial_{x^{\mu}x^{\nu}}^{2}%
\Phi \}_{\mu,\nu=1}^{n}$. Thus $\sum_{k}\mathbb{E}[|\eta_{k}^{N}|]\rightarrow
0$. The rest terms in the summation of the right side of (\ref{Ito}) are
$\xi_{t}^{N}+\zeta_{t}^{N}$, with%
\begin{align*}
\xi_{t}^{N}  &  =\sum_{k=0}^{N-1}\{ \partial_{x^{\mu}}\Phi(X_{t_{k}^{N}%
})[\alpha^{\mu}(t_{k+1}^{N}-t_{k}^{N})+\eta^{\mu}(\left \langle B\right \rangle
_{t_{k+1}^{N}}-\left \langle B\right \rangle _{t_{k}^{N}})+\beta^{\mu
}(B_{t_{k+1}^{N}}-B_{t_{k}^{N}})]\\
&  +\frac{1}{2}\partial_{x^{\mu}x^{\nu}}^{2}\Phi(X_{t_{k}^{N}})\beta^{\mu
}\beta^{\nu}(B_{t_{k+1}^{N}}-B_{t_{k}^{N}})(B_{t_{k+1}^{N}}-B_{t_{k}^{N}})\}
\end{align*}
and
\begin{align*}
\zeta_{t}^{N}  &  =\frac{1}{2}\sum_{k=0}^{N-1}\partial_{x^{\mu}x^{\nu}}%
^{2}\Phi(X_{t_{k}^{N}})[\alpha^{\mu}(t_{k+1}^{N}-t_{k}^{N})+\eta^{\mu
}(\left \langle B\right \rangle _{t_{k+1}^{N}}-\left \langle B\right \rangle
_{t_{k}^{N}})]\\
&  \times \lbrack \alpha^{\nu}(t_{k+1}^{N}-t_{k}^{N})+\eta^{\nu}(\left \langle
B\right \rangle _{t_{k+1}^{N}}-\left \langle B\right \rangle _{t_{k}^{N}})]\\
&  +\beta^{\nu}[\alpha^{\mu}(t_{k+1}^{N}-t_{k}^{N})+\eta^{\mu}(\left \langle
B\right \rangle _{t_{k+1}^{N}}-\left \langle B\right \rangle _{t_{k}^{N}%
})](B_{t_{k+1}^{N}}-B_{t_{k}^{N}}).
\end{align*}
We observe that, for each $u\in \lbrack t_{k}^{N},t_{k+1}^{N})$
\begin{align*}
&  \mathbb{E}[|\partial_{x^{\mu}}\Phi(X_{u})-\sum_{k=0}^{N-1}\partial_{x^{\mu
}}\Phi(X_{t_{k}^{N}})\mathbf{I}_{[t_{k}^{N},t_{k+1}^{N})}(u)|^{2}]\\
&  =\mathbb{E}[|\partial_{x^{\mu}}\Phi(X_{u})-\partial_{x^{\mu}}\Phi
(X_{t_{k}^{N}})|^{2}]\\
&  \leq c^{2}\mathbb{E}[|X_{u}-X_{t_{k}^{N}}|^{2}]\leq C[\delta+\delta^{2}].
\end{align*}
Thus $\sum_{k=0}^{N-1}\partial_{x^{\mu}}\Phi(X_{t_{k}^{N}})\mathbf{I}%
_{[t_{k}^{N},t_{k+1}^{N})}(\cdot)$ tends to $\partial_{x^{\mu}}\Phi(X_{\cdot
})$ in $M_{G}^{2}(0,T)$. Similarly,
\[
\sum_{k=0}^{N-1}\partial_{x^{\mu}x^{\nu}}^{2}\Phi(X_{t_{k}^{N}})\mathbf{I}%
_{[t_{k}^{N},t_{k+1}^{N})}(\cdot)\rightarrow \partial_{x^{\mu}x^{\nu}}^{2}%
\Phi(X_{\cdot})\text{, in \ }M_{G}^{2}(0,T).
\]
Let $N\rightarrow \infty$, by the definitions of the integrations with respect
to $dt$, $dB_{t}$ and $d\left \langle B\right \rangle _{t}$ the limit of
$\xi_{t}^{N}$ in $L_{G}^{2}(\mathcal{F}_{t})$ is just the right hand of
(\ref{B-Ito}). By the estimates of the next remark, we also have $\zeta
_{t}^{N}\rightarrow0$ in $L_{G}^{1}(\mathcal{F}_{t})$. We then have proved
(\ref{B-Ito}).
\end{proof}

\begin{remark}
We have the following estimates: for $\psi^{N}\in M_{G}^{1,0}(0,T)$ such that
$\psi_{t}^{N}=\sum_{k=0}^{N-1}\xi_{t_{k}}^{N}\mathbf{I}_{[t_{k}^{N}%
,t_{k+1}^{N})}(t)$, and $\pi_{T}^{N}=\{0\leq t_{0},\cdots,t_{N}=T\}$ with
$\lim_{N\rightarrow \infty}\mu(\pi_{T}^{N})=0$ and $\sum_{k=0}^{N-1}%
\mathbb{E}[|\xi_{t_{k}}^{N}|](t_{k+1}^{N}-t_{k}^{N})\leq C$, for all
$N=1,2,\cdots$, we have
\[
\mathbb{E}[|\sum_{k=0}^{N-1}\xi_{k}^{N}(t_{k+1}^{N}-t_{k}^{N})^{2}%
|]\rightarrow0,
\]
and, thanks to Lemma \ref{Lem-Qua2},%
\begin{align*}
\mathbb{E}[|\sum_{k=0}^{N-1}\xi_{k}^{N}(\left \langle B\right \rangle
_{t_{k+1}^{N}}-\left \langle B\right \rangle _{t_{k}^{N}})^{2}|]  &  \leq
\sum_{k=0}^{N-1}\mathbb{E[}|\xi_{k}^{N}|\cdot \mathbb{E}[(\left \langle
B\right \rangle _{t_{k+1}^{N}}-\left \langle B\right \rangle _{t_{k}^{N}}%
)^{2}|\mathcal{F}_{t_{k}^{N}}]]\\
&  =\sum_{k=0}^{N-1}\mathbb{E[}|\xi_{k}^{N}|](t_{k+1}^{N}-t_{k}^{N}%
)^{2}\rightarrow0,
\end{align*}
as well as%
\begin{align*}
&  \mathbb{E}[|\sum_{k=0}^{N-1}\xi_{k}^{N}(\left \langle B\right \rangle
_{t_{k+1}^{N}}-\left \langle B\right \rangle _{t_{k}^{N}})(B_{t_{k+1}^{N}%
}-B_{t_{k}^{N}})|]\\
&  \leq \sum_{k=0}^{N-1}\mathbb{E[}|\xi_{k}^{N}|]\mathbb{E[}(\left \langle
B\right \rangle _{t_{k+1}^{N}}-\left \langle B\right \rangle _{t_{k}^{N}%
})|B_{t_{k+1}^{N}}-B_{t_{k}^{N}}|]\\
&  \leq \sum_{k=0}^{N-1}\mathbb{E[}|\xi_{k}^{N}|]\mathbb{E[}(\left \langle
B\right \rangle _{t_{k+1}^{N}}-\left \langle B\right \rangle _{t_{k}^{N}}%
)^{2}]^{1/2}\mathbb{E[}|B_{t_{k+1}^{N}}-B_{t_{k}^{N}}|^{2}]^{1/2}\\
&  =\sum_{k=0}^{N-1}\mathbb{E[}|\xi_{k}^{N}|](t_{k+1}^{N}-t_{k}^{N}%
)^{3/2}\rightarrow0.
\end{align*}
We also have
\begin{align*}
&  \mathbb{E}[|\sum_{k=0}^{N-1}\xi_{k}^{N}(\left \langle B\right \rangle
_{t_{k+1}^{N}}-\left \langle B\right \rangle _{t_{k}^{N}})(t_{k+1}^{N}-t_{k}%
^{N})|]\\
&  \leq \sum_{k=0}^{N-1}\mathbb{E[}|\xi_{k}^{N}|(t_{k+1}^{N}-t_{k}^{N}%
)\cdot \mathbb{E}[(\left \langle B\right \rangle _{t_{k+1}^{N}}-\left \langle
B\right \rangle _{t_{k}^{N}})|\mathcal{F}_{t_{k}^{N}}]]\\
&  =\sum_{k=0}^{N-1}\mathbb{E[}|\xi_{k}^{N}|](t_{k+1}^{N}-t_{k}^{N}%
)^{2}\rightarrow0.
\end{align*}
and
\begin{align*}
\mathbb{E}[|\sum_{k=0}^{N-1}\xi_{k}^{N}(t_{k+1}^{N}-t_{k}^{N})(B_{t_{k+1}^{N}%
}-B_{t_{k}^{N}})|]  &  \leq \sum_{k=0}^{N-1}\mathbb{E[}|\xi_{k}^{N}%
|](t_{k+1}^{N}-t_{k}^{N})\mathbb{E}[|B_{t_{k+1}^{N}}-B_{t_{k}^{N}}|]\\
&  =\sqrt{\frac{2}{\pi}}\sum_{k=0}^{N-1}\mathbb{E[}|\xi_{k}^{N}|](t_{k+1}%
^{N}-t_{k}^{N})^{3/2}\rightarrow0.\
\end{align*}
\end{remark}

We now consider a more general form of It\^{o}'s formula. Consider%
\[
X_{t}^{\nu}=X_{0}^{\nu}+\int_{0}^{t}\alpha_{s}^{\nu}ds+\int_{0}^{t}\eta
_{s}^{\nu}d\left \langle B,B\right \rangle _{s}+\int_{0}^{t}\beta_{s}^{\nu
}dB_{s}.
\]

\begin{proposition}
\label{Prop-Ito}Let $\alpha^{\nu}$, $\beta^{\nu}$ and $\eta^{\nu}$,
$\nu=1,\cdots,n$, are bounded processes of $M_{G}^{2}(0,T)$. Then for each
$t\geq0$ and in $L_{G}^{2}(\mathcal{F}_{t})$ we have%
\begin{align}
\Phi(X_{t})-\Phi(X_{s})  &  =\int_{s}^{t}\partial_{x^{\nu}}\Phi(X_{u}%
)\beta_{u}^{\nu}dB_{u}+\int_{s}^{t}\partial_{x_{\nu}}\Phi(X_{u})\alpha
_{u}^{\nu}du\label{Ito-form1}\\
&  +\int_{s}^{t}[\partial_{x^{\nu}}\Phi(X_{u})\eta_{u}^{\nu}+\frac{1}%
{2}\partial_{x^{\mu}x^{\nu}}^{2}\Phi(X_{u})\beta_{u}^{\mu}\beta_{u}^{\nu
}]d\left \langle B\right \rangle _{u}\nonumber
\end{align}
\end{proposition}

\begin{proof}
We first consider the case where $\alpha$, $\eta$ and $\beta$ are step
processes of the form%
\[
\eta_{t}(\omega)=\sum_{k=0}^{N-1}\xi_{k}(\omega)\mathbf{I}_{[t_{k},t_{k+1}%
)}(t).
\]
From the above Lemma, it is clear that (\ref{Ito-form1}) holds true. Now let
\[
X_{t}^{\nu,N}=X_{0}^{\nu}+\int_{0}^{t}\alpha_{s}^{\nu,N}ds+\int_{0}^{t}%
\eta_{s}^{\nu,N}d\left \langle B\right \rangle _{s}+\int_{0}^{t}\beta_{s}%
^{\nu,N}dB_{s}%
\]
where $\alpha^{N}$, $\eta^{N}$ and $\beta^{N}$ are uniformly bounded step
processes that converge to $\alpha$, $\eta$ and $\beta$ in $M_{G}^{2}(0,T)$ as
$N\rightarrow \infty$. From Lemma \ref{Lem-26}%
\begin{align}
\Phi(X_{t}^{\nu,N})-\Phi(X_{0})  &  =\int_{s}^{t}\partial_{x^{\nu}}\Phi
(X_{u}^{N})\beta_{u}^{\nu,N}dB_{u}+\int_{s}^{t}\partial_{x_{\nu}}\Phi
(X_{u}^{N})\alpha_{u}^{\nu,N}du\label{N-Ito}\\
&  +\int_{s}^{t}[\partial_{x^{\nu}}\Phi(X_{u}^{N})\eta_{u}^{\nu,N}+\frac{1}%
{2}\partial_{x^{\mu}x^{\nu}}^{2}\Phi(X_{u}^{N})\beta_{u}^{\mu,N}\beta_{u}%
^{\nu,N}]d\left \langle B\right \rangle _{u}\nonumber
\end{align}
Since%
\begin{align*}
\mathbb{E[}|X_{t}^{\nu,N}-X_{t}^{\nu}|^{2}]  &  \leq3\mathbb{E[}|\int_{0}%
^{t}(\alpha_{s}^{N}-\alpha_{s})ds|^{2}]+3\mathbb{E[}|\int_{0}^{t}(\eta
_{s}^{\nu,N}-\eta_{s}^{\nu})d\left \langle B\right \rangle _{s}|^{2}]\\
+3\mathbb{E[}|\int_{0}^{t}(\beta_{s}^{\nu,N}-\beta_{s}^{\nu})dB_{s}|^{2}]  &
\leq3\int_{0}^{T}\mathbb{E}[(\alpha_{s}^{\nu,N}-\alpha_{s}^{\nu})^{2}%
]ds+3\int_{0}^{T}\mathbb{E}[|\eta_{s}^{\nu,N}-\eta_{s}^{\nu}|^{2}]ds\\
&  +3\int_{0}^{T}\mathbb{E}[(\beta_{s}^{\nu,N}-\beta_{s}^{\nu})^{2}]ds,
\end{align*}
we then can prove that, in $M_{G}^{2}(0,T)$, we have (\ref{Ito-form1}).
Furthermore
\begin{align*}
\partial_{x^{\nu}}\Phi(X_{\cdot}^{N})\eta_{\cdot}^{\nu,N}+\partial_{x^{\mu
}x^{\nu}}^{2}\Phi(X_{\cdot}^{N})\beta_{\cdot}^{\mu,N}\beta_{\cdot}^{\nu,N}  &
\rightarrow \partial_{x^{\nu}}\Phi(X_{\cdot})\eta_{\cdot}^{\nu}+\partial
_{x^{\mu}x^{\nu}}^{2}\Phi(X_{\cdot})\beta_{\cdot}^{\mu}\beta_{\cdot}^{\nu}\\
\partial_{x_{\nu}}\Phi(X_{\cdot}^{N})\alpha_{\cdot}^{\nu,N}  &  \rightarrow
\partial_{x_{\nu}}\Phi(X_{\cdot})\alpha_{\cdot}^{\nu}\\
\partial_{x^{\nu}}\Phi(X_{\cdot}^{N})\beta_{\cdot}^{\nu,N}  &  \rightarrow
\partial_{x^{\nu}}\Phi(X_{\cdot})\beta_{\cdot}^{\nu}%
\end{align*}
We then can pass limit in both sides of (\ref{N-Ito}) and get (\ref{Ito-form1}).
\end{proof}

\section{Stochastic differential equations}

We consider the following SDE defined on $M_{G}^{2}(0,T;\mathbb{R}^{n})$:
\begin{equation}
X_{t}=X_{0}+\int_{0}^{t}b(X_{s})ds+\int_{0}^{t}h(X_{s})d\left \langle
B\right \rangle _{s}+\int_{0}^{t}\sigma(X_{s})dB_{s},\ t\in \lbrack0,T].
\label{SDE}%
\end{equation}
where the initial condition $X_{0}\in \mathbb{R}^{n}$ is given and
$b,h,\sigma:\mathbb{R}^{n}\mapsto \mathbb{R}^{n}$ are given Lipschitz
functions, i.e., $|\phi(x)-\phi(x^{\prime})|\leq K|x-x^{\prime}|$, for each
$x$, $x^{\prime}\in \mathbb{R}^{n}$, $\phi=b$, $h$ and $\sigma$. Here the
horizon $[0,T]$ can be arbitrarily large. The solution is a process $X\in
M_{G}^{2}(0,T;\mathbb{R}^{n})$ satisfying the above SDE. We first introduce
the following mapping on a fixed interval $[0,T]$:%
\[
\Lambda_{\cdot}(Y):=Y\in M_{G}^{2}(0,T;\mathbb{R}^{n})\longmapsto M_{G}%
^{2}(0,T;\mathbb{R}^{n})\  \
\]
by setting $\Lambda_{t}$ with
\[
\Lambda_{t}(Y)=X_{0}+\int_{0}^{t}b(Y_{s})ds+\int_{0}^{t}h(Y_{s})d\left \langle
B\right \rangle _{s}+\int_{0}^{t}\sigma(Y_{s})dB_{s},\ t\in \lbrack0,T].
\]

We immediately have

\begin{lemma}
For each $Y,Y^{\prime}\in M_{G}^{2}(0,T;\mathbb{R}^{n})$, we have the
following estimate:%
\[
\mathbb{E}[|\Lambda_{t}(Y)-\Lambda_{t}(Y^{\prime})|^{2}]\leq C\int_{0}%
^{t}\mathbb{E}[|Y_{s}-Y_{s}^{\prime}|^{2}]ds,\ t\in \lbrack0,T],
\]
where $C=3K^{2}$.
\end{lemma}

\begin{proof}
This is a direct consequence of the inequalities (\ref{ine-dt}), (\ref{e2})
and (\ref{qua-ine}).
\end{proof}

We now prove that SDE (\ref{SDE}) has a unique solution. By multiplying
$e^{-2Ct}$ on both sides of the above inequality and then integrate them on
$[0,T]$. It follows that%
\begin{align*}
\int_{0}^{T}\mathbb{E}[|\Lambda_{t}(Y)-\Lambda_{t}(Y^{\prime})|^{2}%
]e^{-2Ct}dt  &  \leq C\int_{0}^{T}e^{-2Ct}\int_{0}^{t}\mathbb{E}[|Y_{s}%
-Y_{s}^{\prime}|^{2}]dsdt\\
&  =C\int_{0}^{T}\int_{s}^{T}e^{-2Ct}dt\mathbb{E}[|Y_{s}-Y_{s}^{\prime}%
|^{2}]ds\\
&  =(2C)^{-1}C\int_{0}^{T}(e^{-2Cs}-e^{-2CT})\mathbb{E}[|Y_{s}-Y_{s}^{\prime
}|^{2}]ds.
\end{align*}
We then have
\[
\int_{0}^{T}\mathbb{E}[|\Lambda_{t}(Y)-\Lambda_{t}(Y^{\prime})|^{2}%
]e^{-2Ct}dt\leq \frac{1}{2}\int_{0}^{T}\mathbb{E}[|Y_{t}-Y_{t}^{\prime}%
|^{2}]e^{-2Ct}dt.
\]
We observe that the following two norms are equivalent in $M_{G}%
^{2}(0,T;\mathbb{R}^{n})$:
\[
\int_{0}^{T}\mathbb{E}[|Y_{t}|^{2}]dt\thicksim \int_{0}^{T}\mathbb{E}%
[|Y_{t}|^{2}]e^{-2Ct}dt.
\]
From this estimate we can obtain that $\Lambda(Y)$ is a contract mapping.
Consequently, we have

\begin{theorem}
There exists a unique solution $X\in M_{G}^{2}(0,T;\mathbb{R}^{n})$ of the
stochastic differential equation (\ref{SDE}).
\end{theorem}

\section{Appendix}

For $r>0$, $1<p,q<\infty$ with $\frac{1}{p}+\frac{1}{q}=1$, we have
\begin{align}
|a+b|^{r}  &  \leq \max \{1,2^{r-1}\}(|a|^{r}+|b|^{r}),\  \  \forall
a,b\in \mathbb{R}\label{ee4.3}\\
|ab|  &  \leq \frac{|a|^{p}}{p}+\frac{|b|^{q}}{q}. \label{ee4.4}%
\end{align}

\begin{proposition}%
\begin{align}
\mathbb{E}[|X+Y|^{r}]  &  \leq C_{r}(\mathbb{E}[|X|^{r}]+\mathbb{E[}%
|Y|^{r}]),\label{ee4.5}\\
\mathbb{E}[|XY|]  &  \leq \mathbb{E}[|X|^{p}]^{1/p}\cdot \mathbb{E}%
[|Y|^{q}]^{1/q},\label{ee4.6}\\
\mathbb{E}[|X+Y|^{p}]^{1/p}  &  \leq \mathbb{E}[|X|^{p}]^{1/p}+\mathbb{E}%
[|Y|^{p}]^{1/p}. \label{ee4.7}%
\end{align}
In particular, for $1\leq p<p^{\prime}$, we have $\mathbb{E}[|X|^{p}%
]^{1/p}\leq \mathbb{E}[|X|^{p^{\prime}}]^{1/p^{\prime}}.$
\end{proposition}

\begin{proof}
(\ref{ee4.5}) follows from (\ref{ee4.3}). We set
\[
\xi=\frac{X}{\mathbb{E}[|X|^{p}]^{1/p}},\  \  \eta=\frac{Y}{\mathbb{E}%
[|Y|^{q}]^{1/q}}.
\]
By (\ref{ee4.4}) we have%
\begin{align*}
\mathbb{E}[|\xi \eta|]  &  \leq \mathbb{E}[\frac{|\xi|^{p}}{p}+\frac{|\eta|^{q}%
}{q}]\leq \mathbb{E}[\frac{|\xi|^{p}}{p}]+\mathbb{E}[\frac{|\eta|^{q}}{q}]\\
&  =\frac{1}{p}+\frac{1}{q}=1.
\end{align*}
Thus (\ref{ee4.6}) follows. We now prove (\ref{ee4.7}):
\begin{align*}
\mathbb{E}[|X+Y|^{p}]  &  =\mathbb{E}[|X+Y|\cdot|X+Y|^{p-1}]\\
&  \leq \mathbb{E}[|X|\cdot|X+Y|^{p-1}]+\mathbb{E}[|Y|\cdot|X+Y|^{p-1}]\\
&  \leq \mathbb{E}[|X|^{p}]^{1/p}\cdot \mathbb{E[}|X+Y|^{(p-1)q}]^{1/q}\\
&  +\mathbb{E}[|Y|^{p}]^{1/p}\cdot \mathbb{E[}|X+Y|^{(p-1)q}]^{1/q}%
\end{align*}
We observe that $(p-1)q=p$. Thus we have (\ref{ee4.7}).
\end{proof}

\end{document}